\documentclass[reqno]{gokova}

\begin{document}

\newtheorem{thm}{Theorem}[section]
\newtheorem{cor}[thm]{Corollary}
\newtheorem{lem}[thm]{Lemma}
\newtheorem{prop}[thm]{Proposition}
\newtheorem{ax}{Axiom}
\newtheorem{conj}[thm]{Conjecture}

\theoremstyle{definition}
\newtheorem{defn}{Definition}[section]

\theoremstyle{remark}
\newtheorem{rem}{\rm\bfseries{Remark}}[section]
\newtheorem*{notation}{Notation}

\newtheorem{ques}{\rm\bfseries{Question}}[section]
\newtheorem{cons}[rem]{\rm\bfseries{Construction}}
\newtheorem{exm}[rem]{\rm\bfseries{Example}}

%\numberwithin{equation}{section}

\title{Evidence for a conjecture of Pandharipande}

\author{Jim Bryan}
\thanks{The author is supported by an Alfred P. Sloan Research
Fellowship and NSF grant DMS-0072492}

\address{Department of Mathematics, Tulane University}
\email{jbryan@math.tulane.edu}

%\newtheorem{thm}{Theorem}[section]
%\newtheorem{theorem}[thm]{Theorem}
%\newtheorem{conjecture}[thm]{Conjecture}

%==========================================
%% Do not edit the following command
\volume{7}
%==========================================

\maketitle 

\section{Introduction}

In \cite{Pa}, Pandharipande studied the relationship between the
enumerative geometry of certain 3-folds and the Gromov-Witten invariants.
In some good cases, enumerative invariants (which are manifestly integers)
can be expressed as a rational combination of Gromov-Witten
invariants. Pandharipande speculated that the same combination of
invariants should yield integers even when they do not have any enumerative
significance on the 3-fold. In the case when the 3-fold is the product of a
complex surface and an elliptic curve, Pandharipande has computed this
combination of invariants on the 3-fold in terms of the Gromov-Witten
invariants of the surface \cite{Pa-comm}. This computation yields surprising
conjectural predictions about the genus 0 and genus 1 Gromov-Witten
invariants of complex surfaces. The conjecture states that certain rational
combinations of the genus 0 and genus 1 Gromov-Witten invariants are always
integers. Since the Gromov-Witten invariants for surfaces are often
enumerative (as oppose to 3-folds), this conjecture can often also be
interpreted as giving certain congruence relations among the various
enumerative invariants of a surface.

In this note, we state Pandharipande's conjecture and we prove it for an
infinite series of classes in the case of $\mathbf{CP}^{2}$ blown-up at 9
points. In this case, we find generating functions for the numbers
appearing in the conjecture in terms of quasi-modular forms
(Theorem~\ref{thm: computation of A and B}). We then prove the integrality
of the numbers by proving a certain a congruence property of modular forms
that is reminiscent of Ramanujan's mod 5 congruences of the partition
function (Theorem~\ref{thm: mod 10 congruence for G}).

\section{The conjecture}

Let $X$ be a smooth complex projective surface (or more generally, a
symplectic 4-manifold), let $K$ be its canonical class, and let $\chi (X)$
be its Euler characteristic. Let $\beta\in H_{2} (X,\mathbf{Z})$ and let $g
(\beta)$ be defined by $2g (\beta)-2=\beta\cdot (K+\beta)$. Define $c
(\beta)$ to be $-\beta \cdot K$ and assume that $c (\beta)>0$. Let $N^{r}
(\beta)$ be the genus $r$ Gromov-Witten invariant of $X$ in the class
$\beta$ where we have imposed $c (\beta)+r-1$ point constraints. By
convention we will say $N^{r} (0)=0$.

\begin{conj}[Pandharipande]
Define $a (\beta )$ by 
\[
a (\beta )=-\frac{1}{12}g (\beta )N^{0} (\beta )
\]
and 
define $b (\beta )$ by 
\begin{align*}
b (\beta  )=&\quad \frac{1}{2880} \left(12g (\beta)^{2}+g (\beta)c (\beta)-24g (\beta) \right)N^{0} (\beta)\\
&+\frac{1}{240}\chi (X)N^{1} (\beta)\\
&+\frac{1}{240}\sum _{\beta'+\beta''=\beta} \binom{c (\beta)-1}{c (\beta')} (\beta'\cdot \beta'')
(\beta''\cdot \beta '')N^{1} (\beta')N^{0}(\beta'') .
\end{align*}
Then  $a (\beta )$ and $b (\beta )$ are integers.
\end{conj}

\begin{rem}
This conjecture is related to the proposal of Gopakumar and Vafa that
relates the Gromov-Witten invariants of Calabi-Yau 3-folds to conjecturally
integer valued invariants (``BPS state counts'', or ``BPS
invariants''). Pandharipande has generalized the Gopakumar-Vafa formula to
Fano classes in non-Calabi-Yau 3-folds (see \cite{Pa}). In this
formulation, the numbers $a (\beta )$ and $b (\beta )$ are respectively
genus 1 and genus 2 ``BPS invariants'' for the surface cross an elliptic
curve. The reason that these are expressible in terms of ordinary
Gromov-Witten invariants of the surface is that the Hodge class in
$\overline{M}_{g}$ (which appears in the computation of the virtual class)
is readily expressible in terms of boundary classes for $g=1$ and
$g=2$. For arbitrary $g$ there will also be predictions for the invariants
of the surface, but they will involve gravitational descendants in general.
\end{rem}

\section{The case of $\mathbf{CP}^{2}$ blown-up at 9 points}

Let $X$ be $\mathbf{CP}^{2}$ blown up at nine points. Let $F=-K$ be the
anti-canonical class and let $S$ be the exceptional divisor of one of the
blow-ups (so if $X$ is elliptically fibered, then $F$ is the fiber and $S$
is a section). Let $\beta_{n}=S+nF$. Then $N^{r} (\beta_{n})$ was computed
in \cite{BL}. We will find a nice generating functions for the numbers $a
(\beta_{n})$ and $b (\beta _{b})$ and will prove that they are integers
thus verifying Pandharipande's conjecture for $X$ for this infinite series
of classes.

Note that $c (\beta_{n})=1$, $g (\beta_{n})=n$, and $\chi (X)=12$. Since
for $N^{0} (\beta'')$ to be non-zero, we need $c (\beta'')=1$, the sum must
have $c (\beta'')=1$ and $c (\beta')=0$. It follows that $\beta''$ and
$\beta'$ are of the form $S+kF$ and $(n-k)F$ respectively. Thus we have
\begin{align*}
a (\beta _{n})=&-\frac{1}{12}nN^{0} (\beta _{n})\\
b (\beta_{n})=&\quad \frac{1}{2880} (12n^{2}-23n)N^{0} (\beta_{n})\\
&+\frac{1}{20}N^{1} (\beta_{n})\\
&+\frac{1}{240}\sum _{k=0}^{n-1} (n-k) (2k-1)N^{1} ((n-k)F)N^{0} (\beta_{k}).
\end{align*}

Define 
\begin{align*}
A (q)&=\sum _{n=0}^{\infty }a (\beta _{n})q^{n},\\
B (q)&=\sum _{n=0}^{\infty }b (\beta _{n})q^{n}.
\end{align*}
We will find an expression for $A (q)$ and $B (q)$ in terms of
quasi-modular forms. Let $\sigma (k)=\sum _{d|k}d$ and let $p (k)$ be the
number of partitions of $k$. Define
\begin{align*}
G (q)&=\sum _{k=1}^{\infty }\sigma (k)q^{k},\\
P (q)&=\sum _{k=1}^{\infty }p (k)q^{k}\\
&=\prod _{m=1}^{\infty } (1-q^{m})^{-1},\\
P_{\alpha } (q)&= (P (q))^{\alpha },\\
D&=q\frac{d}{dq}.
\end{align*}
Note that $G$ and $P$ are closely related to well known (quasi-) modular forms:
$G-1/24$ is the Eisenstein series $G_{2} $ and $q^{1/24}P_{-1}$ is the Dedekind
$\eta $ function.

%Define $D=q\frac{d}{dq}$ and $P_{\alpha } (q)= (P (q))^{\alpha }$. 
With
this notation, the results of \cite{BL} (Theorem 1.2) give
\begin{align*}
\sum _{n=0}^{\infty }N^{0} (\beta_{n})q^{n} &=P_{12}\\
\sum _{n=1}^{\infty }N^{1} (\beta_{n})q^{n}&=P_{12}DG.
\end{align*}

Furthermore, one can show that 
\[
N^{1} (lF)=\frac{1}{l}\sigma (l)
\]
(when the blow-up points are generic, this comes from the multiple covers
of the unique elliptic curve in the class $F$). We thus have
\begin{align*}
A (q)&=-\frac{1}{12}DP_{12}\\
B (q)&= \frac{1}{2880} (12D^{2}-23D)P_{12} + \frac{1}{20}P_{12}DG\\
& +\frac{1}{240}\sum _{n\geq 1}\sum _{k=0}^{n-1} (2k-1)\sigma
(n-k)N^{0} (\beta_{k})q^{n-k}q^{k}\\
&= \frac{1}{240}D^{2}P_{12}-\frac{23}{2880}DP_{12}+\frac{1}{20}P_{12}DG\\
& +\frac{1}{240}\sum _{m\geq 1}\sum _{k\geq 0} (2k-1)\sigma (m)N^{0}
(\beta_{k})q^{k}q^{m}\\
&=\frac{1}{240}D^{2}P_{12}-\frac{23}{2880}DP_{12}+\frac{1}{20}P_{12}DG
+\frac{1}{240}G (2DP_{12}-P_{12})
\end{align*}

Now, by a standard calculation, $G=P_{-1}DP$ and so
$DP_{12}=12P_{12}G$. Substituting and simplifying we arrive at:
\begin{thm}\label{thm: computation of A and B}
The following equations holds:
\begin{align*}
A (q)&=-P_{12}\cdot G\\
B (q)&=\frac{1}{10}P_{12}\left\{7G^{2}-G +DG\right\}.
\end{align*}
\end{thm}

This theorem immediately shows that the coefficients of $A$ are
integers. On the other hand, the integrality of the coefficients of $B$
requires the following theorem:
\begin{thm}\label{thm: mod 10 congruence for G}
The following equation holds:
\[
7G^{2}-G+DG\equiv 0 \pmod {10}.
\]
\end{thm}
\textsc{Proof:}
By a simple calculation mod 5, we have:
\[
7G^{2}-G+DG\equiv 3P_{-2} (D^{2}-D)P_{2} \pmod{5}
\]
and so to prove that the above expression is 0 mod 5, it suffices to prove
that $(D^{2}-D)P_{2}\equiv 0 \pmod{5}$. Using the Jacobi triple product
formula and the Euler inversion formula, it is easy to show that the $k$th
coefficient of $P_{2}=P_{-3}P_{5}$ is divisible by 5 unless $k$ is 0 or 1
mod 5 (see \cite{Drost}). In other words:
\[
P_{2} (q)\equiv r (q^{5})+qs (q^{5}) \pmod 5.
\]
It follows that $DP_{2}\equiv qs (q^{5})\pmod 5$ and so $D^{2}P_{2}\equiv
DP_{2}\pmod 5$ as desired.

On the other hand, it is easy to compute that 
\[
7G^{2}-G+DG\equiv P_{-1} (D^{2}+D)P\pmod 2.
\]
This expression is 0 mod 2 since the $k$th coefficient of $(D^{2}+D)P$ is
$k(k+1)p (k)$.

Thus we have established that $7G^{2}-G+DG$ is 0 mod 2 and mod 5 and so the
theorem is proved. \qed

\end{document}